\documentclass[12pt,a4paper]{amsart}
\usepackage[english]{babel}
\usepackage[utf8]{inputenc}
\usepackage[T1]{fontenc}

\usepackage{textcase}

\usepackage{pdflscape}
\usepackage{amsmath}
\usepackage{amssymb}
\usepackage{latexsym}
\usepackage{graphics}
\usepackage{epsfig}
\usepackage{color}

\usepackage[colorlinks=true, allcolors=blue]{hyperref}

\usepackage[a4paper,hmargin=2.4cm,vmargin=2.5cm]{geometry}

\newtheorem{theorem}{Theorem}
\newtheorem*{notmytheorem}{Theorem}
\newtheorem*{lemma*}{Lemma}
\newtheorem{prop}{Proposition}

\newtheorem{remark}{Remark}

\newtheorem*{oc*}{Open question}

\DeclareMathOperator{\Ncl}{Ncl}

\title[On the word problem for just infinite groups]{On the word problem \\ for just infinite groups }
\author{Alexey Talambutsa}
\address{\parbox{\linewidth}{Steklov Mathematical Institute of RAS, \\ Gubkina str., 8, 119991 Moscow, Russia\\
HSE University, Laboratory of Theoretical Computer Science, \\ Pokrovsky Boulevard, 11, 109028 Moscow, Russia}}
\email[]{altal@mi-ras.ru}


\begin{document}

\begin{abstract}
We establish several results on the word problem for just infinite groups. 
First, for finitely generated just infinite groups we show that the word problem is uniformly decidable for presentations with recursively enumerable sets of relations. Our proof does not use the Wilson--Grigorchuk theorem on the classification of just infinite groups and proceeds directly from the definition, using ideas from classical results on decidability of the word problem: Kuznetsov's theorem and McKinsey--Maltsev theorem.

For countably generated presentations of just infinite groups with a recursively enumerable set of relations we show that the word problem is decidable in all cases except locally finite groups without a computable lower bound on subgroup sizes.  %The decidability is also established in the case when the group is locally finite, but there exists a computable lower bound on the sizes of the subgroups $\langle a_1,\ldots,a_n\rangle$, where $S=\{a_1,a_2,\ldots\}$ is the countable generating set of the presentation.

Finally, we construct presentations of countably generated locally finite groups with recursively enumerable set of relations, for which the word problem is undecidable. Yet, there exist other presentations of these groups, for which the word problem is decidable.

\end{abstract}

\maketitle

\section{Introduction}

A group $G$ is called \emph{just infinite} if it is infinite, and for every proper normal subgroup $H\triangleleft G$ the quotient group $G/H$ is finite. The class of just infinite groups was introduced by J.~Wilson in the early 1970s (see \cite{Wilson1971}). Clearly, this class contains all infinite simple groups and the infinite cyclic and dihedral groups. According to the notable result of J. Mennicke (see \cite{Mennicke1965}), the projective matrix groups $PSL(n,\mathbb{Z})$ possess this property for $n\geq 3$. In the 1980s R.~Grigorchuk and Gupta--Sidki constructed many examples of just infinite groups acting on rooted trees (see \cite{Grigorchuk1980, Grigorchuk1985} and \cite{GuptaSidki1983}). Periodic hereditarily just infinite groups (i.e., residually finite groups for which every subgroup of finite index is just infinite) were found in the 2010s by J.~Wilson in \cite{Wilson2010} and by R.~Grigorchuk and P.~Shumyatskii in \cite{GriShum2017}.

\smallskip

In Grigorchuk's paper \cite{Grigorchuk2000}, the class of branch groups acting on rooted trees was defined, and based on \cite{Wilson1971}, the following result, often called the trichotomy theorem, was proved:

\begin{notmytheorem}
\label{classification}
Let $G$ be a finitely generated just infinite group. Then either $G$ is a branch group, or $G$ contains a normal subgroup of finite index isomorphic to a direct product of finitely many copies of a group $L$, where $L$ is either simple or hereditarily just infinite.
\end{notmytheorem}

Since the decidability of the word problem for a finitely generated group and for its subgroup of finite index are equivalent, it follows from Theorem \ref{classification} that, in order to solve the word problem in just infinite groups, it suffices to solve it for three subclasses: branch, simple, and residually finite groups that are hereditarily just infinite.

The algorithmic decidability of the word problem for finitely generated simple groups with a recursively enumerable set of relations was established by A.V.~Kuznetsov (see \cite{Kuznetsov1958}; cf.~\cite[Theorem 3.6]{LyndonSchupp}). A positive solution to the word problem for residually finite groups was obtained by A.I.~Maltsev \cite{Maltsev1958} who used an earlier method of J.C.C.~McKinsey \cite{McKinsey1943}. Later, the same result was also obtained by V.~Dyson and A.~Mostowski (see \cite{Dyson1964} and \cite{Mostowski1966}). At the same time, as was shown by V.~Dyson and S.~Meskin (see \cite{Dyson1974} and \cite{Meskin1974}), there exist examples of groups from this class with recursively enumerable sets of relations for which this problem is algorithmically undecidable.

An algorithm solving the word problem for spinal groups --- which include the Grigorchuk and Gupta–Sidki groups, as well as many other branch groups --- was constructed in \cite{BarGrigSun2003}. However, this algorithm takes as input not a presentation of the group in terms of generators and relations, but a special oracle that is related to the action of generating automorphisms of the tree. We also note that for the series of Grigorchuk groups $G_{\omega}$ it was established that the word problem in $G_{\omega}$ is decidable if and only if the sequence $\omega$ is computable. Since it is not known whether every finitely generated branch group is spinal (in particular, this is unknown for the group in \cite{Neumann1986}), a complete answer to the question of decidability of the word problem for the class of branch groups has not been established.

\smallskip

The goal of this paper is to investigate the decidability of the word problem for finitely and countably generated just infinite groups with a recursively enumerable set of relations. We will consider the \emph{word problem}, when the presentation of the group is fixed, and the \emph{uniform word problem}, for which it is provided as part of the algorithm input.

\smallskip

For a finitely generated group $G$, a \emph{recursively enumerable presentation} is a pair $(S,\mathcal{R})$ consisting of some finite alphabet of generators $S=\{a_1,a_2,\ldots,a_k\}$, and an algorithm $\mathcal R$ that enumerates the defining relations of the group $G$ in the generators $S$. We will also use the notation $\langle S \mid \mathcal{R} \rangle$ for the presentation of the group $G$, assuming that in this case $\mathcal R$ describes a set of relations. In the case when the group is generated by a countably infinite set $S=\{a_1,a_2,\ldots\}$, the same notation as above will be used; it is clear that with a suitable encoding the presentation $(S,\mathcal{R})$ can be encoded as a finite word in a finite alphabet.

\section{Finitely generated just infinite groups}

The main result of this section is the following statement:

\begin{theorem}
\label{fg_theorem}
    There exists an algorithm that, given a recursively enumerable presentation $\langle S \mid \mathcal{R} \rangle$ of a finitely generated just infinite group $G$ and a word $x$ in the generators $S$, decides whether the equality $x=1$ holds in the group $G$.
\end{theorem}

\begin{proof}
We will use the same scheme that is used both in Kuznetsov's theorem and in the McKinsey--Maltsev theorem. Two algorithms work in parallel, one of which tries to prove the equality $x=1$ in the group $G$, and the other tries to refute this equality. Since exactly one of the two algorithms terminates, a definitive answer is obtained in finite time.

\smallskip

\emph{Algorithm 1} (to prove $x=1$). The reasoning here is similar to both Kuznetsov's theorem and the McKinsey–Maltsev theorem. Suppose that the algorithm $\mathcal{R}$ has enumerated, up to the present moment, the defining relations $R_1, R_2, \ldots, R_l$. Then we have $l$ parallel processes $P_1,P_2,\ldots,P_l$, and the process $P_s$ will be enumerating all freely reduced words that are equal to $1$ in the group $\langle S \mid R_1=1, R_2=1, \ldots, R_s=1 \rangle$. These will be Dyck words that have the form
\begin{equation}
\prod_{j=1}^s t_j R^{\varepsilon_j}_{i_j} t_j^{-1}=1, \quad \text{where } t_j\in F_k \text{ and } \varepsilon_{j}\in \{-1,0,1\}.
\label{dyck_word}
\end{equation}
Let us fix an ordering for the complete enumeration of the exponents $\varepsilon_{j}$ and the reduced words $t_j$ from the free group $F_k$ of rank $k$. It is known that every word that is equal to $1$ in the group $G=\langle S \mid \mathcal{R} \rangle$, after performing free group cancellations, coincides graphically with some Dyck word (see \cite[Lemma 5.1]{AdianDurnev2000}). If $x=1$ then let $q$ be the largest index $i_j$ that was used for $x$ in the product \eqref{dyck_word}. Then the word $x$ will be found in any enumerating process $P_s$ with $s\geq q$, and once $x$ is found, the whole algorithm will terminate. If $x\ne 1$, then Algorithm 1 will never terminate.

\medskip

\emph{Algorithm 2} (to prove $x\ne 1$). Here we apply reasoning that combines two different ideas that are used in Kuznetsov's theorem with the McKinsey--Maltsev theorem.

First, as in Kuznetsov's theorem, we extend the presentation of the group $G$ with the relation $x=1$ and obtain a new group $G_1$ with a recursively enumerable presentation $G_1=\langle S \mid \{x\} \cup \mathcal R\rangle$. Note that if $x=1$ in $G$, then $G_1 \cong G$ is infinite. Conversely, if $x \neq 1$, the just infinite property implies that the quotient $G_1 = G/\Ncl(x)$ is finite.

We now proceed, as in the McKinsey--Maltsev theorem, to enumerate all possible multiplication tables for finite groups $T_1,T_2,\ldots$ and try to find an epimorphism $\psi:T_i\to G_1$. Clearly, such an epimorphism will exist if and only if $G_1$ is a finite group. To check the existence of an epimorphism, we consider each table $T_q$ as a presentation of a group, taking all its elements $u_1,u_2,\ldots,u_r$ as generators and all equalities of the form $u_iu_j=u_k$ from the multiplication table as relations. In order to find the elements that map to the generators $S$ under the epimorphism, we will enumerate all possible maps $\tau:\{u_1,u_2,\ldots,u_r\}\to S^{\pm *}$ such that the image of $\{u_1,u_2,\ldots,u_r\}$ includes $\{a_1,a_2,\ldots,a_k\}$. For each such map $\tau$ we will try to derive from the relations $\{x\}\cup \mathcal R$ of the group $G_1$ the relations
\begin{equation}
\label{checking_equations}
\tau(u_i)\tau(u_j)=\tau(u_k), \, \text{where } u_iu_j=u_k \text{ is an equality from the multiplication table } T_q.
\end{equation} Since, as mentioned above, the set of all relations of $G_1$ is enumerable, if the relations \eqref{checking_equations} hold in $G_1$, then all of them will be derived eventually 
%\footnote{Note that in the Dyson--Mostowski theorem, a suitable homomorphism $\psi:G \to T_i$ is sought in the reverse direction, and then the correctness of $\psi$ cannot be checked in finite time for an infinite recursively enumerable presentation of the group $G$. In particular, this is why, unlike in Kuznetsov's theorem, the presentation of the group is assumed to be finite.}
. Then for the given multiplication table $T_q$ and substitution of generators $\tau$ we obtain that the group $G_1$ is a quotient group of the group $T_q$, and therefore it is finite. On the other hand, if $G_1$ is finite, then its generators $a_1,a_2,\ldots,a_k$ partition $S$ into disjoint sets of equal elements, and a suitable surjection $\tau$ from the generators of the multiplication table of the group $G_1$ onto the set $S$ will yield the required epimorphism. It remains to note that if the group $G_1$ is infinite, then no suitable epimorphism exists, and Algorithm 2 does not terminate.
\end{proof}

Boone and Higman asked (see \cite{BooneHigman1974}; cf.\cite{BBMZ2025}) whether a finitely generated group $G$ has a decidable word problem if and only if $G$ can be embedded into a finitely presented simple group. One direction of it is Kuznetsov's theorem, and the other is still open. In view of Theorem~\ref{fg_theorem} one may suggest the following weakened version of this question: 
\begin{oc*}
\normalfont{
Does a finitely generated group $G$ have a decidable word problem if and only if $G$ can be embedded into finitely presented just infinite group?
}
\end{oc*}
However, the main difficulty which arises on the way to an answer seems to us the same as for the Boone--Higman conjecture: it is hard both to construct finitely presented infinite simple groups and finitely presented just infinite groups. In particular, there are no known examples of finitely presented branch groups.

\section{Countably generated groups}
\label{sec-cgg}

\subsection{Uniform word problem}
We start by showing that Theorem~\ref{fg_theorem} cannot be generalized to the case of a countable generating set. In this case the general problem becomes undecidable even if one considers only the presentations of the infinite cyclic group.

\begin{theorem}
\label{cg_g_theorem}
    There is no algorithm that, given a countably generated recursively enumerable presentation $G=\langle S \mid \mathcal{R} \rangle$ of an infinite cyclic group and a word $x$ in the generators $S$, decides whether the equality $x=1$ holds in the group $G$.
\end{theorem}
\begin{proof}
Consider a universal Turing machine $UTM(M,w)$ taking as input a Turing machine $M$ and a binary input word $w$. Let $e(M)$ be some effective encoding function that maps finite programs of Turing machines into binary words.

Now, let $M$ be a Turing machine. We will associate to $M$ an infinite recursively enumerable set of relations $\mathcal R_M$ in the generating set $S=\{a_1,a_2,a_3\ldots\}$. To do this, with the use of $UTM$ we execute the machine $M$ on the input word $e(M)$ step by step. If the machine $M$ has not stopped after the current step $k$, we write out the relation $a_{k+1}=1$. If at some step $d$ the machine $M$ stops, we complement $\mathcal R_M$ with the relation $a_1=1$ and infinitely continue the recursive enumeration with $a_{d+2}=1, a_{d+3}=1,\ldots$. Clearly, a program for computing $M$ can be transformed into the program which enumerates $\mathcal R_M$.

\smallskip

If $M$ does not halt on input $e(M)$, the described procedure will enumerate the relations $\{a_2=1,a_3=1,\ldots\}$. If the machine stops after some step $d$, we have that the only generator of $S$ which is not equal to one is $a_{d+1}$, whilst in the non-halting case it is $a_1$. Thus, deciding whether the relation $a_1=1$ holds in the group $\langle S \mid \mathcal R\rangle$ is equivalent to checking if the Turing machine is self-terminating, which is known to be undecidable. It is clear that in any case, the corresponding group is infinite cyclic.
\end{proof}

\begin{remark}\normalfont{
One can easily generalize the statement of Theorem~\ref{cg_g_theorem} to any non-trivial finitely presented group. In particular, for a cyclic group of order two one can make it by adding extra relations $a^2_1=1$ and $a^2_d=1$ to the presentation. Thus, the uniform word problem is also undecidable for the countably generated recursively presented simple groups.
}
\end{remark}

\subsection{Word problem for a fixed presentation} Even though, as shown above, the uniform word problem is undecidable in the case of countably infinite generation, we can prove that the word problem for a fixed presentation of a group is decidable for some classes of these groups. First, we observe that Kuznetsov's theorem can be generalized (see also \cite[Proposition 3.2]{BBMZ2025}) .

\begin{theorem}
\label{cg_s_theorem}
    Let $G=\langle S \mid \mathcal{R} \rangle$ be a countably generated recursively enumerable presentation of a simple group. There exists an algorithm that given a word $x$ in the generators $S$, decides whether the equality $x=1$ holds in the group $G$.
\end{theorem}
\begin{proof}
We need to check, whether $x=1$. As in the standard proof of Kuznetsov's theorem, we first note that the set of all words which are equal to $1$ is recursively enumerable. Thus, if $x=1$, we will obtain $x$ at some moment of this enumeration. In the second parallel process we can factorize $G$ by an additional relation $x=1$ and check if the group $G/\Ncl(x)$ is trivial. The standard proof of Kuznetsov's theorem does that by checking whether all generators become trivial, but in our case the set $S$ is infinite, so we would not be able to finish such a test in finite time. However, it is sufficient to make a test with just one element. Let us choose a word $w$ such that $w\ne 1$ in $G$
%\footnote{We have no algorithmic way of finding this word $w$, which makes the proof non-constructive.}
; now, instead of considering the generators, we can check the equality $w=1$ in the factor group $G/\Ncl(x)$. If $x\ne 1$, then $G/\Ncl(x)$ is a trivial group, so $w=1$ in it, and we will eventually see $w$ in the enumeration of the words that are equal to $1$ in $G/\Ncl(x)$.
\end{proof}

We will see later that the analogous result is not true for the class of all just infinite groups, but at first we show that it is valid for some subclasses.

\begin{theorem}
\label{cg_jinlf_theorem}
    Let $G=\langle S \mid \mathcal{R} \rangle$ be a countably generated recursively enumerable presentation of a just infinite group, which is \textbf{not} locally finite. There exists an algorithm that, given a word $x$ in the generators $S$, decides whether the equality $x=1$ holds in $G$.
\end{theorem}
\begin{proof}
Since $G$ is not locally finite, there exists $d$ such that the subgroup $H=\langle a_1,a_2,\ldots,a_d\rangle$ is infinite. Again, as in Theorem~\ref{cg_s_theorem}, it is sufficient to detect $x\ne 1$. Due to the just infinite property for $G$, the latter can be done once we find an algorithm that stops if and only if the image $\varphi(H)$ in the factorization $\varphi: G\to G/\Ncl(x)$ is finite.

Since $\varphi(H)$ is a finitely generated subgroup of a recursively presented group, it is itself recursively presented with a finite number of generators. This is well-known for finitely generated ambient groups (see \cite[Lemma 2.1]{Miller}), but can be easily generalized to the case of countable generation. In our case it is straightforward, as we can consider $G/\Ncl(x)$ as a factor group of the infinitely generated free group $F_S=\langle S \rangle$ by the normal closure of recursive set of words $\mathcal R\cup\{x\}$, and then obtain $\varphi(H)$ as the quotient group
\[
\varphi(H)=\frac{\langle a_1,\ldots,a_d\rangle_{F_S}}{\langle a_1,\ldots,a_d\rangle_{F_S}\cap\Ncl_{F_S}(\mathcal R\cup\{x\})},
\]
so one can obtain the recursive presentation for $\varphi(H)$ by intersecting the enumeration of trivial words in the group $G/\Ncl(x)$ with the set of all words in the alphabet $\{a_1,a_2,\ldots,a_d\}^{\pm}$.

Now, having the recursive presentation of the group $\varphi(H)$ we apply the same argument as in Algorithm 2 of Theorem~\ref{fg_theorem} and list all finite groups and their surjective mappings onto $\varphi(H)$. If $\varphi(H)$ is finite, one multiplication table would suffice to stop the algorithm, and so the check $x \ne 1$ would be finished.
\end{proof}

In the case when the group is locally finite, there are still subclasses, when the word problem with fixed presentation of a group is decidable.

\begin{theorem}
\label{cg_jilf_theorem}
    Let $G=\langle a_1,a_2,\ldots \mid \mathcal{R} \rangle$ be a countably generated recursively enumerable presentation of a locally finite just infinite group. Then the word problem for $G$  is decidable if and only if there exists a computable non-decreasing unbounded function $f:\mathbb N\to \mathbb N$ such that $|\langle a_1,\ldots,a_k\rangle_G|\geq f(k)$.
\end{theorem}
\begin{proof}
If the word problem is decidable, then for any $k$ one can enumerate the ball $B_k(r)$ of radius $r$ (i.e. the elements of length up to $r$) in the finitely generated subgroup $\langle a_1,\ldots,a_k\rangle_G$. Since this subgroup is finite, the balls in it eventually coincide with the subgroup; thus $|B_k(R(k))|=|B_k(R(k)+1)|=|\langle a_1,\ldots,a_k\rangle_G|$ for some $R(k)\in \mathbb Z_{\geq 0}$. The function $|B_k(R(k))|$
is computable by the above description, and thus serves as the required function $f(k)$.

\smallskip

Now suppose the function $f(k)$ exists, and we need to show that the word problem is decidable. As in Theorem~\ref{cg_jinlf_theorem}, we will be checking the inequality $x\ne 1$ by factoring the group $G$ by the normal closure of element $x$. If $x\ne 1$, then the group $G/\Ncl(x)$ is finite, and its subgroups $H_k=\langle a_1, a_2,\ldots,a_k\rangle_{G/\Ncl(x)}$ are finite too. Then for repeatedly increasing parameter $k=1,2,\ldots$ we start running in parallel a separate procedure, which enumerates the presentations of $H_k$ the same way we did in Theorem~\ref{cg_jinlf_theorem}.

At each stage $s$ the enumerated part of the presentation of $H_k$ defines a finitely presented group $Q_{k,s}$, which is an extension of $H_k$. Even though we are not guaranteed that the group $Q_{k,s}$ is finite, we run a process that searches for a finite group $T_k$ and an epimorphism $T_k\to Q_{k,s}$.  Once such a group $T_k$ is found, we can check whether $|T_k|<f(k)$; if this happens, the inequality
\begin{equation}
|\langle a_1,\ldots,a_k\rangle_G|\geq f(k)>|T_k|
\label{size_est}
\end{equation}
forces $T_k$ to be the proper quotient of the subgroup $\langle a_1,\ldots,a_k\rangle_G$. The latter shows that the factorization $G\to G/\Ncl(x)$ is non-trivial, hence $x\ne 1$.

\smallskip

Now we need to show that if $x\ne 1$, the inequality \eqref{size_est} holds true for some group $Q_{k,s}$. Since the group $G/\Ncl(x)$ is finite, we have $|G/\Ncl(x)|=m$ for some $m\in \mathbb N$. Since the function $f(k)$ is unbounded, then for some $k_0\in \mathbb N$ one has
\[
f(k_0)>m
.\]

Clearly, $|H_{k_0}|\leq m$ since $H_{k_0}$ is a subgroup in $G/\Ncl(x)$. Due to Maltsev compactness theorem for the presentation of the finite group $H_{k_0}$, there exists a moment $s$ such that the partial presentation $Q_{k_0,s}$ coincides with the group $H_{k_0}$. This means that $|Q_{k_0,s}|=|H_{k_0}|<m$. By listing all finite groups which possess epimorphisms onto the group $Q_{k_0,s}$, eventually we will find the multiplication table $T_{k_0}$ which defines a finite group isomorphic to $Q_{k_0,s}$, so then $|T_{k_0}|=|Q_{k_0,s}|=|H_{k_0}|<m$. Thus, for $k=k_0$ the inequality \eqref{size_est} is true.
\end{proof}

\begin{remark}
\label{alg_finite}
\normalfont{
One can prove that the required function $f(k)$ in Theorem~\ref{cg_jilf_theorem} exists, once there exists an algorithm that enumerates an infinite set of distinct elements. Indeed, let $f(k)$ be the number of group words written in the first $k$ generators among the first $k$ words produced by such an algorithm. Clearly, this function is computable and non-decreasing. Since the algorithm $A$ enumerates the infinite set, the function $f(k)$ will not be bounded. Due to the definition of $f$, the inequality $f(k)\leq|\langle a_1,\ldots,a_k\rangle_G|$ holds true. (Note that non-existence of such an algorithm for finitely generated groups is called \emph{algorithmic finiteness}, and its presence makes the solution of the word problem extremely complicated (see \cite{MyasnikovOsin})).
%This argument implies that any countably generated recursively presented group $\langle S \mid \mathcal{R}\rangle$ has another recursive presentation, in which is has decidable word problem. Indeed, let $S=\{a_1,a_2,\ldots \}$ and choose an increasing function $t(n)$ which enumerates the subset $V\subset S$ of the elements which are
}
\end{remark}

\subsection{Monolithic groups}

The classification of just infinite groups in \cite{Wilson1971} goes through a natural strict dichotomy between the residually finite groups and monolithic groups, that is established as follows. Recall that for any group $G$ its \emph{monolith} $M(G)$ is the intersection of all non-trivial normal subgroups in $G$. The group $G$ is called \emph{monolithic} if $M(G)\ne \{1\}$. Since any non-trivial normal subgroup $H$ in a just infinite group $G$ has finite index, for such a group $G$ one has
\[M(G)=\bigcap \{ \,H \mid H \triangleleft G, \ [G:H]<\infty\},\]
and the latter intersection is trivial if and only if the group $G$ is residually finite.

\smallskip

To the best of the author's knowledge, the following result has not appeared in the literature. Clearly, it can be seen as a generalization of Theorem \ref{cg_s_theorem}. 

\begin{theorem}
\label{cg_m_theorem}
    Let $G=\langle S \mid \mathcal{R} \rangle$ be a countably generated recursively enumerable presentation of a monolithic group. Then there exists an algorithm that, given a word $x$ in the generators $S$, decides whether the equality $x=1$ holds in $G$.
\end{theorem}
\begin{proof}
Again, as before, we only need to detect if $x\ne 1$. Choose a non-trivial element $w\in M(G)$. It holds that $x\ne 1$ if and only if $\Ncl(x)$ is a non-trivial (normal) subgroup in $G$, which in turn happens if and only if $w\in \Ncl(x)$. The subgroup $\Ncl(x)<G$ is a countable union of finitely generated subgroups $\langle x^{y_1},\ldots,x^{y_s}\rangle_G$, where all finite subsets $\{y_1,\ldots,y_s\}\subset G$ are considered. Since the set of words in the alphabet $S^{\pm}$ that are equal to $1$ in $G$ is recursively enumerable, then so is the set $Z$ of words in the alphabet $S^{\pm}$ representing $\Ncl(x)$. Thus, to detect $x\ne 1$ we can run the enumeration of the set $Z$ and check if we find the word $w$ in it.
\end{proof}

\section{Countable presentations with undecidable word problem}
The goal of this section is to show that there exist presentations of just infinite groups with countably infinite set of generators and recursively enumerable set of relations, for which the word problem is undecidable. The results obtained in Section~\ref{sec-cgg} show that any just infinite group admitting such a presentation must be locally finite and residually finite.
%First, we will show that the residual finiteness is necessary for constructing just infinite examples with undecidability.
%Theorem~\ref{cg_jinlf_theorem} and Theorem~\ref{cg_m_theorem} show that if there is a recursively enumerable presentation $G=\langle S \mid \mathcal R \rangle$ of a just infinite group with undecidable word problem, then the group has to be locally finite and residually finite. 

%Moreover, from Remark~\ref{alg_finite} one can see that for such a presentation the indices of generators which are equal to $1$ is a \emph{simple set}, which, contrary to its name, is a rather complicated object. Recall that a set $C\subset \mathbb N$ is a \emph{simple set} if it is recursively enumerable and has an infinite complement $\overline{C}$ that does not contain infinite recursively enumerable subsets\footnote{The interested reader can see that Theorem~\ref{cg_jilf_theorem} shows that, moreover, this set of indices is hypersimple.}. Below we construct the presentations of just infinite groups satisfying all these restrictions.

\subsection{Shaded sets}
Our first goal is to construct a subset of natural numbers that is recursively enumerable, not recursive, and has a specific origin, which we call a shaded set.

\smallskip

If $a,b\in \mathbb N$ with $a\leq b$, we write $[a,b]$ for the subset $\{ a, \ldots, b \}\subset \mathbb N$, and call it a \emph{segment}.

\smallskip

Let $\alpha$ be a recursively enumerable sequence of segments $\alpha=\{[l_t,r_t]\}_{t=0}^\infty$, for which we construct a \emph{shaded set} $\mathcal I(\alpha)=\bigcup_{t=0}^{\infty} I_t$, where $I_t$ is a segment defined as follows:
\begin{equation}
\label{it-case}
I_t =
 \begin{cases}
   [l_t,r_t] &\text{if $l_t>\max \{r_i \mid i<t\}$},\\
   [l_t,\max \{r_i \mid i\leq t\}] &\text{otherwise}.
 \end{cases}
\end{equation}
This construction can be described as follows: The first condition means that the segment $[l_t,r_t]$ is to the right of all previous segments in $\alpha$, and in this case, we include it into the union $\mathcal I(\alpha)$; otherwise the segment $[l_t, r_t]$ is absorbed into a larger interval extending to the maximal previously reached right border, which can be thought of as casting a ``big shade''.

\smallskip

We will also use the same definition $\mathcal I(\alpha)$ for a finite sequence of segments $\alpha=\{[l_t,r_t]\}_{t=0}^n$.

\begin{prop}
\label{prop-shaded-set}
There exists a recursively enumerable sequence of segments $\alpha=\{[l_t,r_t]\}_{t=0}^\infty$ for which the shaded set $\mathcal I(\alpha)$ is recursively enumerable and not recursive.
\end{prop}
\begin{proof}
To construct the corresponding sequence $\alpha$ we will be using the Busy Beaver function $BB(n):\mathbb N\to \mathbb N$, which equals the maximum number of non-blank symbols an $n$-state Turing machine can print on an initially blank tape before eventually halting. It is known that $BB(n)$ grows faster than any computable function, i.e., if the function $f:\mathbb N\to \mathbb N$ possesses the property $f(n)\geq BB(n)$ for any $n$, then $f$ is not computable.

To define the sequence $\alpha$ we will run following algorithm, which computes a function that estimates $BB$ from below: At major step $k$ one runs the first $k$ steps of all Turing machines containing at most $k$ states, which gives us numbers $LBB(1;k), LBB(2;k), \ldots, LBB(k;k)$ equal to the maximum number of non-blank symbols that were written by the machines that stopped within $k$ steps of execution. We will assume that $LBB(n,k)=0$ for $n>k$, and thus have a total function $LBB:\mathbb N\times \mathbb N\to \mathbb N$. Clearly, the function $LBB$ is computable, non-decreasing in the second argument and subject to the inequality $LBB(n,k)\leq BB(n)$ that is true for any $n,k\in \mathbb N$. By construction of $LBB$ one also has
\begin{equation}
\label{lim_inequality}
\lim_{k\to\infty} LBB(n,k) = BB(n).
\end{equation}

Now we define the sequence of segments $\alpha$. We define the segment $[l_0,r_0]$ to be $[1,1]$ and continue defining further segments along the major steps of the algorithm which computes the values of the function $LBB$, which we described above. We will be appending the sequence $\alpha=\{[l_t,r_t]\}_{t=0}^\infty$ such that after each major step $k$ the segments indexed from $0$ to some $M_k\in \mathbb N$ are defined, and the shaded set of this fragment $\alpha_k=\{[l_t,r_t]\}_{t=0}^{M_k}$ can be expressed as a union of increasing segments having form:
\begin{equation}
\label{invariant}
\mathcal I(\{[l_t,r_t]\}_{t=0}^{M_k})=[1,v_1]\cup[1+v_1,v_2]\cup\ldots \cup [1+v_{k-1},v_k],
\end{equation}
and their size is bounded by the inequality
\begin{equation}
\label{inv_inequality}
    v_i-v_{i-1}\geq LBB(i,k) \quad \text{for} \quad 1\leq i\leq k.
\end{equation}

At each major step $k+1$ we look at the length of all segments in \eqref{invariant} and compare the size $p_i=v_i-v_{i-1}$ of the $i$-th segment with the newly computed number $LBB(k+1,i)$. If $p_i\geq LBB(k+1,i)$, we do nothing. If $p_i<LBB(k+1,i)$, we stop the comparison process and start modifying our invariant \eqref{invariant} by appending segments to $\alpha$ in the following way:

First, we add to $\alpha$ a segment $[1+v_{i-1},1+W]$, where $W=\max(v_k,LBB(k+1,i))$. Clearly, this new segment absorbs all segments in \eqref{invariant} numbered from $i$ to $k$, so after this modification the shaded set of the extended finite part of $\alpha$ contains exactly $i$ segments. After that, for each $j\in \{i+1,\ldots, k\}$ we append to $\alpha$ a new segment of length $\max(p_{j},LBB(k+1,j))$, and we separate these $k-i$ segments inside $\mathbb N$ by singleton gaps. Finally, we also append to the sequence $\alpha$ a segment of length $LBB(k+1,k+1)$, which becomes the segment numbered $k+1$ in the decomposition \eqref{invariant} for the extended sequence $\alpha$ of length $M_{k+1}$.
%\footnote{Since at step $k+1$ we append to $\alpha$ not more than $k+1$ new segments, one has $M_k\leq (k+1)k/2$.}

\smallskip

The described process of appending segments to the finite sequence $\alpha_k$ at the step $k+1$ by its construction preserves the structure \eqref{invariant} and the inequality \eqref{inv_inequality}. The sequence $\alpha$ which is obtained as a limit is recursively enumerable, and hence so is the shaded set $\mathcal{I}(\alpha)$, which we are now going to analyze.

Due to the equality \eqref{lim_inequality} and inequality \eqref{inv_inequality}, at some step $k$ the length of the first interval in \eqref{invariant} stops changing after its length equals or exceeds $BB(1)$, and no further overwriting of the first segment will be happening. By the same argument, after this, at some moment the second segment in \eqref{invariant} stabilizes after its length reaches $BB(2)$, then the third, etc. It follows, that for some increasing sequence of natural numbers $1=u_0<u_1<u_2<\ldots$ we have
\begin{equation}
\label{lim-str}
\mathcal I(\alpha)=\bigcup_{i=0}^{\infty}[1+u_i,u_{i+1}] \quad \text{with } u_i-u_{i-1}\geq BB(i) \text{ for any } i\in \mathbb N.
\end{equation}
From this we establish that the set $\mathcal{I}(\alpha)$ is not recursive. Otherwise, one would be able to check natural numbers in the increasing order for their presence in $\mathcal{I}(\alpha)$ and allocate the $i$-th segment in \eqref{lim-str} by finding the singletons gaps. The latter allows to find the value of the function $D(i)=u_{i}-u_{i-1}$ for any $i$. Since the inequality $D(i)\geq BB(i)$ holds for any $i\in \mathbb N$, the function $D(i)$ is not computable, so we get to a contradiction.
\end{proof}

\subsection{Shaded presentations of the iterated wreath products}
The groups $W_C$ which we consider are particular iterated wreath products in the series, described by J.~Wilson in \cite{Wilson1971} as examples of residually finite locally finite just infinite groups. 

\smallskip

Let $C$ be a non-trivial perfect finite group and let $\mathbb Z_-$ denote the set of negative integers. The infinite wreath product $W_C=\text{Wr}_{i\in \mathbb Z_-}C$ in the sense of P.~Hall (see \cite{Hall1962}) can also be seen as a direct limit of permutational wreath products $W_C=\lim_{n\to \infty} G_n$, where $G_1=C$ and $G_{n+1}=G_n \wr C$. As shown in \cite[Construction 2]{Wilson1971}, the group $W_C$ is just infinite. We note that if the permutational group $C$ acts on a set of cardinality $q$, then the group $W_C \cong W_C \wr C$ is an infinitely generated self-similar branch group acting on a regular tree of valence $q$.

\begin{theorem}
For each group $W_C$ there exists a countably generated recursively enumerable presentation $\langle S \mid \mathcal U \rangle$ with an undecidable word problem.
\label{th-undecidable}
\end{theorem}
\begin{proof}
First, we describe the standard way to construct for the group $W_C$ a recursively enumerable presentation $\langle S \mid \mathcal R \rangle$ with decidable word problem. Then we will show how to modify this process using the recursively enumerable sequence $\alpha$ from Proposition~\ref{prop-shaded-set} to obtain the desired group presentation $\langle S \mid \mathcal U \rangle$ with undecidable word problem.

\smallskip

Let us choose a finite presentation $C=\langle a_1,\ldots,a_d \mid q_1=1, \ldots,q_m=1 \rangle$ for the initial perfect finite permutational group $C$ acting on the set of size $r$. 

Due to the associativity of permutational wreath products, one can as well write the group $W_C$ as $\lim_{n\to \infty} G_n$, where $G_0=\{1\}$ and $G_{n+1}= C \wr G_{n}$. Thus, $G_{n+1}$ is a semidirect product $B_n \rtimes G_n$, where the normal subgroup $B_n$ is a direct product of $|q^n|$ copies of the group $C$. The permutational group $G_n$ acts on $B_n$ by permuting these copies. 

To obtain a presentation for the group $G_{n+1}$ we take a presentation of $G_n$ and extend it to a presentation of the free product $G_{n}*B_n$ by adding $dq^n$ new generators and relations in them: the $mq^n$ defining relations for $q^n$ copies of the group $C$ and commutations between the generators of different copies. After that we add finitely many relations that define the action of generators for $G_n$ on the copies of $C$ and thus obtain $G_{n+1}$ as a quotient of $G_{n}*B_n$. In this construction the canonical projecting homomorphism $\pi_{n+1}: G_{n+1}=B_n \rtimes G_n \to G_n$ is realized by sending the $dq^n$ new generators to $1$. 

If we define the recursive function
\begin{equation}
f(0)=0, \ f(1)=d \quad \text{and} \quad  f(n+1)=f(n)+dq^n,
\end{equation}
then the process above gives a presentation of the group $G_n$ with $f(n)$ generators, and we will denote it as $GP_n=\langle s_1,s_2,\ldots,s_{f(n)} \mid R_n \rangle$. Note that the presentation $GP_n$ can be considered as the initial part of the presentation $GP_{n+1}$, which also agrees with the canonical embedding $\psi_n:G_n\to G_{n+1}$. Passing to the limit as $n\to\infty$, one obtains a countable set of generators and a countable set of relations that define the limit presentation $W_C=\langle S \mid \mathcal R \rangle$.

%The description of the presentation $\langle S \mid \mathcal R \rangle$ above shows that $G_n\cong\langle s_1,\ldots, s_{E(n)} \mid R_n \rangle$, where $R_n$ are some first relations of $\mathcal R$ that are written in the alphabet $\{ s_1,\ldots, s_{E(n)}\}^\pm$.

%\smallskip

%Consider $k,n\in \mathbb N$ such that $k<n$. A presentation of the group $G_k$ can be obtained from the presentation $GP_n=\langle s_1,\ldots, s_{f(n)} \mid R_n \rangle$ by adding relations $s_{E(k)+1}=1,\ldots,s_{E(n)}=1$, which corresponds to a composition of projections $\pi_n\circ\ldots\circ\pi_{k+1}=:\pi_{(n,k)}$. %This also allows us to see $W_C$ as the inverse limit of groups $G_n$ %We will be using quotients of this type and extensions of the presentations $G_n\to G_{n+1}$ described above to recursively define the presentation $\langle S \mid \mathcal U \rangle$. 

\smallskip

Now let $\alpha=\{[l_t,r_t]\}_{t=0}^\infty$ be a recursively enumerable sequence of segments from Proposition~\ref{prop-shaded-set}. We will use $\alpha$ to build the required presentation $\langle S \mid \mathcal U \rangle$ as a certain limit
%\footnote{As we consider both embeddings and quotients in the construction, the correct limit is inverse, but to make it work, one should extend generators of each presentation $P_k$ to a countable set $S=\{s_1,s_2\ldots \}$.} 
of an infinite sequence of finite presentations $\{P\}_{i=0}^{\infty}$. Each presentation $P_k$ is written in formal generators $\{s_1,s_2,\ldots,s_{f(c(k))}\}$, where $c(k)=\max\{r_i \mid i\leq k\}.$ The presentation $P_k$ presents a wreath product $G_{w(k)}$ for some function $w(k): \mathbb Z_{\geq 0} \to \mathbb Z_{\geq 0}$, which follows to be recursive from our construction, and subject to the inequality $w(k)\leq c(k)$. The group $G_{w(k+1)}$ is obtained from $G_{w(k)}$ by the extension composition $\psi_{w(k)}\circ\ldots\circ\psi_{w(k+1)-1}$ when $w(k+1)\geq w(k)$ or by the projection composition $\pi_{w(k)}\circ\ldots\circ\pi_{w(k+1)+1}$ when $w(k+1)<w(k)$. Before we proceed to the details, we introduce a partition of natural numbers into segments: 
\[
\mathbb N = \bigcup_{i=1}^{\infty} S_i, \ \ \text{where } S_i=[f(i-1)+1,f(i)].
\]
We will use this partition to split generators of $W_C=\langle \mathcal S \mid \mathcal U \rangle$ so that the generator indices of each base group in the iterated wreath product $W_C=\text{Wr}_{i\in \mathbb Z_-}C$ fit into a single segment.

%is the following correspondence between the generators of $P_k$ that are equal to $1$ and the shaded set of sequence $\alpha$:
%\begin{equation}
%\{ i \mid s_i=1 \text{ in } P_k \} = \bigcup \bigl\{ [E(p)+1,E(p+1)] \mid  p\in\mathcal{I}(\{[l_t,r_t]\}_{t=0}^k)  \bigr\}.
%\end{equation}

\smallskip

We start with a trivial group $G_0$ with the empty presentation $P_0$. At each step $k\in \mathbb Z_{\geq 0}$ we consider the segment $[l_k,r_k]$ and similarly to the cases in \eqref{it-case} we define $P_{k+1}$ using $P_k$. 

The key property we preserve is the following: $P_k$ has $f(c(k))=f(\max \{r_i \mid i\leq k\})$ generators, and for each segment $S_i$ either $s_j=1$ in $P_k$ for all $j\in S_i$ when $i\in \mathcal I(\{[l_t,r_t]\}_{t=0}^k)$; or some initial generators $s_j$ with $j\in S_i$ are not equal to $1$ in $P_k$, when $i\notin \mathcal I(\{[l_t,r_t]\}_{t=0}^k)$. 

\emph{Case 1: $l_{k+1}>\max \{r_i \mid i<k+1\}=c(k)$.} In this case we extend the presentation $P_k=\langle s_1,\ldots,s_{c(k)} \mid Q_k \rangle$ by adding new generators and relations. As $P_k$ presents the group $G_{w(k)}$, we extend this group so that $P_{k+1}$ presents $G_{w(k+1)}$, where $w(k+1)=w(k)+(l_{k+1}-c(k)-1)$. To do this, for each $i\in \{1,\ldots,l_{k+1}-c(k)-1\}$ we use the extension map $\psi_{w(k)+i}$ considering new generators $s_t$, where $t\in S_{c(k)+i}$. The first $dq^{w(k)+i-1}$ elements in this segment will be used for indexing the generators of the base group $B_{w(k)+i}$, and map the generators $s_t$ with remaining indices in $S_{w(k)+i}$ to identity. Since $|S_{c(k)+i}|=f(c(k)+i)-f(c(k)+i-1)=dq^{c(k)+i-1}$, and $w(k)\leq c(k)$, this allocation is possible. We also add the relations for the corresponding extension of $G_{w(k)+i-1}$ to $G_{w(k)+i}$ in the standard way. Finally, for each $i\in [l_{k+1},r_{k+1}]$ we add the relations $\{ s_j=1 \mid j \in  S_{i} \}$ ensuring  $c(k+1)=r_{k+1}$. It is easy to see that the key property is preserved by our construction and that the presentation $P_{k+1}$ defines the group $G_{w(k+1)}$.

\emph{Case 2: $l_{k+1}\leq\max \{r_i \mid i<k+1\}$.} Here we only need to factor out some generators in the presentation $P_k$ so that the key property is preserved. To do this, we consider generators $s_j$ with $j\in \{S_i \mid i\in [l_{k+1},c(k)]\}$ and add to $P_{k+1}$ a relation $s_j=1$ for each of them. Since the generators of the top base groups are fit into segments, the described relations factor out some top base groups in $G_{w(k)}$ seen as the iterated wreath product. We have that $c(k+1)=c(k)$ and the key property is preserved by \eqref{it-case} and our construction. One finds the number $w(k+1)$ by finding the size of the set $[1,c(k+1)]\setminus \mathcal I(\{[l_t,r_t]\}_{t=0}^{k+1})$.

\smallskip
Now we can take the countable set of generators $S=\{s_1,s_2,\ldots \}$ and consider the presentation $\langle S \mid \mathcal U\rangle$, where $\mathcal U = \bigcup_{k=1}^{\infty} Q_k$. Since the set $\mathbb N\setminus \mathcal{I}(\alpha)$ is infinite, there will be infinitely many $k\in \mathbb N$ such that the presentation $P_k=\langle s_1,s_2,\ldots,s_{f(c(k))} \mid Q_k \rangle$ remains stable in the union, so we have $G_k<\langle S \mid \mathcal U \rangle$. By considering an increasing subsequence of such indices $k(t)$ we obtain $\langle S \mid \mathcal U \rangle = \lim_{t\to \infty} \langle S\mid Q_{k(t)} \rangle \cong W_C$. 

It remains to notice that the key property of the construction implies that the word problem of $\langle S \mid \mathcal U\rangle$ is reducible to the presence problem for the set $\mathcal{I}(\alpha)$. Indeed, if $i$ is the first element in the segment $S_j$, then $s_i=1$ is equivalent to checking whether $j\in\mathcal I(\alpha)$. Since the latter is undecidable by the choice of the sequence $\alpha$, the result follows.
\end{proof}

\begin{remark}\normalfont{Let us take $C=\langle a_1,\ldots,a_d \mid q_1=1, \ldots,q_m=1 \rangle$ so that all generators are distinct group elements. Then
the presentation $W_C=\langle S \mid \mathcal U \rangle$ from Theorem \ref{th-undecidable} has the following interesting property: the indices of generators which are equal to $1$ in $W_C$ form a \emph{Post's~simple~set}, that is defined as a recursively enumerable set having an infinite complement without infinite recursively enumerable subsets. The second property follows from the observation that all non-trivial elements in the generating set $S$ are distinct, so according to Remark~\ref{alg_finite}, one cannot enumerate any its subset due to the undecidability of the word problem for $\langle S \mid \mathcal U \rangle$.
%\footnote{The interested reader can see that Theorem~\ref{cg_jilf_theorem} shows that, moreover, this set of indices is hypersimple.}
}
\end{remark}

%Moreover, from Remark~\ref{alg_finite} .

%Let $\mathcal E$ be the set of all total computable functions $\mathbb N \to \mathbb N$. We introduce the following operation $\tau: \mathcal E \to \mathcal E$. If $f\in \mathcal E$ and

\section{Acknowledgments}
The author is grateful to Lev D. Beklemishev, Rostislav I. Grigorchuk and John Wilson for useful discussions. The author is indebted to Jordi Delgado, Paul-Henry Leemann and Enric Ventura for their valuable comments on an earlier shorter version of this paper. The work was funded within the framework of the HSE University Basic Research Program.
%The author's work was performed

\end{document}